\documentclass{amsproc}
\usepackage{amsmath,amssymb,amsthm}
\usepackage{graphics}

\newtheorem{thm}{Theorem}[section]

\newtheorem{lem}[thm]{Lemma}

\theoremstyle{definition}

\theoremstyle{remark}

\numberwithin{equation}{section}
\usepackage{enumerate}
\usepackage{graphicx}

\begin{document}

\title{On Triples of Numbers}
\author{Christoph Hering}
\address{
Institute of Mathematics, University of T\"ubingen, 
Auf der Morgenstelle 10, 72076 T\"ubingen, Germany}
\email{hering@uni-tuebingen.de}


\begin{abstract} Among three natural numbers there is always one which is larger than or equal to the Nim sum of the 
remaining two numbers. This amazing fact has many applications.
\end{abstract}



\maketitle

\section{Introduction}\label{Introduction}
On the set of natural numbers we have the usual order, the usual addition and also the somewhat exotic Nim addition. In this paper we investigate the relationship between the usual ordering and the curious Nim addition.

Let $a$ and $b$ be natural numbers (i.e. non-negative integers). Let $a = (a_0 , a_1 , \cdots )$ and $b = (b_0 , b_1, \cdots ) $  be the respective binary expansions $($so, e.g.,$~ a = a_0 + a_1 2^1 + a_2 2^2 + \cdots ).$ Then we define
\[
a\oplus b = a_0\oplus b_0 + (a_1\oplus b_1) 2^1 + (a_2 \oplus b_2) 2^2 + \cdots ,
\]
where $0 \oplus 0 = 1\oplus 1 = 0 $ and $ 0\oplus 1=1\oplus 0 = 1.~a\oplus b$ is called the \emph{Nim sum} of the natural numbers $a$ and $b$. Note that in this paper $0$ is included in the set of natural numbers. The set of natural numbers will be denoted by $\mathbb{N}$.

\section{Triangles of numbers}\label{Sect.3}

Let $a,b$ and $c$ be natural numbers. We call a vertex of the triangle $(a,b,c) $ $ \emph{large} $ if it is larger than the Nim sum of the remaining two numbers. Thus, e.g., a is large if and only if 
\[
a>b \oplus c.
\]
In the same way $a$ is called
\begin{enumerate}[]
  \item $ \emph{aligned} $ if and only if $a=b \oplus c$ , and 

  \item $ \emph{small} $ if and only if $a<b \oplus c$ . 
\end{enumerate}
Correspondingly for the other vertexes. 

If $a$ is aligned, than $a=b \oplus c$, $b=a \oplus c$ and $c=a \oplus b$, and hence all vertexes of the triangle $(a,b,c)$ are aligned. In this case we call the triangle $(a,b,c)$ $ \emph{flat} $ .

\begin{lem} \label{lem:1}Let $a , b $ and $c$ be natural numbers. If $a$ is aligned with respect to the triangle $(a,b,c)$ , then $b$ and $c$ are aligned also. 
\end{lem}

From now on we assume that $a \neq b \oplus c$, i.e. the triangle $(a,b,c)$ is not flat. Let $a = (a_0 , a_1 , \cdots ),~ b = (b_0 , b_1, \cdots ) $ and $~c = (c_0 , c_1 , \cdots )$ be the respective binary expansions $($so, e.g.,$~ a = a_0 + a_1 2^1 + a_2 2^2 + \cdots ). $ Let $j$ be the largest number such that $ a_j \neq b_j \oplus c_j$. 

\vspace{0,2cm}
We have the following possibilities 
\[
\begin{matrix}
a_j	&b_j&c_j	& &a		&b		&c		\\
1		&1	&1		& &large&large&large\\
1		&1	&0		&*\\
1		&0	&1		&*\\
1		&0	&0		& &large&small&small\\
0		&1	&1		&*\\
0		&1	&0		& &small&large&small\\
0		&0	&1		& &small&small&large\\
0		&0	&0		&*\\
\end{matrix}
\]
Here $*$ stands for contradiction. In these cases $a_j=b_j \oplus c_j$, contradicting the definition of $j$.

(We use the fact that for $i>j$ we have $a_i = b_i \oplus c_i$ and therefore $b_i = a_i \oplus c_i$ and $c_i = a_i \oplus b_i$. So, e.g., $b$ and $a \oplus c$ have equal coefficients to the right of the $j$th coefficient.)

Note that the triangle $(a,b,c)$, which is not flat by assumption, always contains a large vertex. In fact, we have got 

\begin{lem} \label{lem:2} Let $a,b$ and $c$ be natural numbers such that $c \neq a \oplus b$ . Then the number of large vertexes of the triangle $(a,b,c)$ is 1 or 3. 

\end{lem}

We call the triangle $(a,b,c)$ \emph{tight} if the number of its large vertexes is 3 and \emph{loose} if this number is only 1. 
We can rephrase Lemma \ref{lem:2}.

\begin{lem} \label{lem:3} Among three natural numbers there is one which is larger than the Nim sum of the remaining two numbers, unless each of the given numbers is equal to the Nim sum of the respective remaining two. 
\end{lem}
We rephrase once more.
\begin{lem} \label{lem:4}
Let $a_1, a_2, a_3 \in \mathbb{N}.$ We can renumber the $a_i$ in such a way that $a_1 \geq a_2 \oplus a_3$.
\end{lem}
Remark. Let $a, b \in \mathbb{N}$ and 
\[
\mathcal{X} = \{\alpha \oplus b \mid \alpha \in \mathbb{N} ~ and ~ \alpha <a\} \cup \{a \oplus \beta \mid \beta \in \mathbb{N} ~ and ~ \beta < b\}. 
\]
Then $a \oplus b \notin \mathcal{X} $ because $(\mathbb{N}, \oplus )$ is a group. Thus $a \oplus b \in \mathbb{N} \setminus \mathcal{X},$ and by Lemma \ref{lem:2}, actually $a \oplus b $ is the smallest element of $\mathbb{N} \setminus \mathcal{X}.$ (Let $c<a \oplus b $. Then the triangle $(a,b,c)$ is not flat and hence contains a large vertex, say $a$ (see Lemma \ref{lem:2}). This implies that $b \oplus c < a$ and hence $c= (b \oplus c) \oplus b \in \mathcal{X}$ and $c\notin \mathbb{N}\setminus \mathcal{X}$.) 

Hence we can say, cum grano salis, that the Nim addition $\oplus $ is the "smallest" binary operation on $\mathbb{N} $ 
leading to a group. In fact it is the "smallest" binary operation on $\mathbb{N} $ with unique solvability of equations (i.e. $a\oplus x = a\oplus y $ implies $x=y$, and $x \oplus a = y \oplus a $ implies $x = y $). 

This been proved first and in a more general situation by Conway in his fabulous book \cite{Co}.

The Three Number Lemma \ref{lem:2} can be used to solve various further first choice construction problems.
 These will be described in a forthcoming paper \cite{KHE}.

\end{document}